\documentclass{article}
\usepackage{ijcai17}

\usepackage{times}

\usepackage{amssymb}
\usepackage{mathtools}
\usepackage{amsmath}
\usepackage{gensymb}
\usepackage{algorithm}
\usepackage{algorithmic}
\usepackage{graphicx}
\usepackage{graphicx}
\usepackage[utf8]{inputenc}
\usepackage{listings}
\usepackage{color}
\usepackage{upquote}
\usepackage{titlesec}
\usepackage{lipsum}
\usepackage{float}
\usepackage{mathtools}
\usepackage{caption}
\usepackage{subcaption}
\usepackage{amsmath}
\usepackage{relsize}
\usepackage{graphicx}
\usepackage[utf8]{inputenc}
\usepackage{listings}
\usepackage{color}
\usepackage{upquote}
\usepackage{titlesec}
\usepackage{lipsum}
\usepackage{float}
\usepackage{mathtools}
\usepackage{caption}
\usepackage{subcaption}

\definecolor{dkgreen}{rgb}{0,0.6,0}
\definecolor{gray}{rgb}{0.5,0.5,0.5}
\definecolor{mauve}{rgb}{0.58,0,0.82}

\title{A comparison between the Shooting and Finite-Difference Method in solving a Nonlinear Boundary Value Problem found in the context of light propagation}
\author{Luke Taylor\\ 
University of Cape Town\\
Department of Mathematics and Applied Mathematics\\
Cape Town, South Africa\\
tylchr011@uct.ac.za}

\begin{document}

\maketitle

\begin{abstract}
The shooting and finite-difference method are both numeric methods that approximate the solution of a BVP to a given accuracy. In this report both methods were implemented in Matlab and compared to each other on a BVP found in the context of light propagation in nonlinear dielectrics. It was observed that the finite-difference method is numerically more stable and converges faster than the shooting method.
\end{abstract}

\section{Introduction}
A Boundary Value Problem \textit{(BVP)} \cite{numericAnalysis} is a differential equation with initial conditions specified at the extremes on the domain of the independent variable. In practise BVPs are solved using numeric techniques as they are accurate (solutions can be found given a defined error threshold) and \textit{fast} (compared to solving complicated BVPs analytically).\\
\begin{equation}
\label{eq:1}
\centering \displaystyle \frac{d^2v}{dr^2}+\frac{1}{r}\frac{dv}{dr}-v+2v^3=0 \quad \frac{dv}{dr}(0)=0, \quad \lim_{r\to\infty} v(r)=0
\end{equation}\\
BVP \ref{eq:1} originally appeared in the context of light propagation in nonlinear dielectrics \cite{lightprop}; Specifically in self-focusing of light beams in the Kerr medium. In this report this differential equation is solved using two popular numeric methods: The Shooting method \cite{shooting} and the Finite-Difference \cite{finitedifference} method for nonlinear problems. Two solutions were sough after, being the monotonically decaying and the one-node solution.

\section{Numeric Approximation Methods}
\begin{figure}[t]
    \centering
    \includegraphics[width=0.5\textwidth]{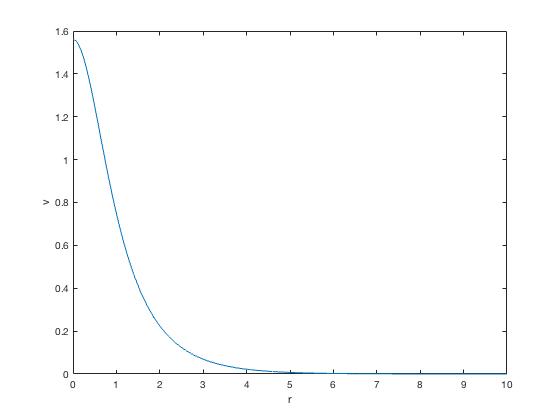}
    \caption{Decaying Solution}
    \label{fig:decay}
\end{figure}
\begin{figure}[t]
    \centering
    \includegraphics[width=0.5\textwidth]{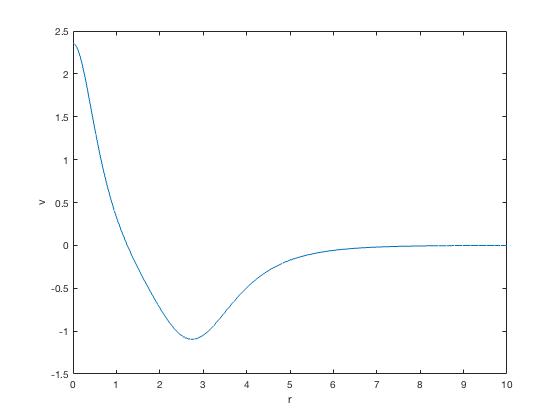}
    \caption{One-Node Solution}
    \label{fig:onenode}
\end{figure}
Numeric Approximations are deployed when an analytical solution is too time consuming to derive and the approximation error is acceptable or when the analytical solution is impossible to derive. The Shooting and Finite-Difference Method are both numeric methods that can approximate second-order boundary-value problems, \\
\begin{equation}
\label{eq:2}
\centering \displaystyle y''=p(x)y'+q(x)y+r(x), \quad \textrm{for} \quad a \leq x \leq b
\end{equation}\\
that are linear or nonlinear with two boundary conditions. Our BVP \ref{eq:1} is nonlinear and hence their iterative nonlinear versions were used. Both algorithmic implementations were completed using Matlab where the source code can be found in the Appendix section.\\
For both methods the solutions were found over the interval $[0, 10]$ (see Discussion for reason) with a max iterations limit of $100$ (the error threshold varied for the different experiments).

\subsection{Nonlinear Shooting Method}
For this method a series of Initial Value Problems (IVPs) are solved using a parameter \textit{p}. The BVP \ref{eq:1} is converted to an IVP by applying the following conditions $v(0)=p$ and $\frac{dv}{dr}(0)=0$. Parameter \textit{p} is iteratively chosen such that $\lim_{k\to\infty} v(\infty, p_k)=v(\infty)=0$. \\
For this implementation \textit{p} was iteratively chosen using the bisection method. The decaying solution (Figure \ref{fig:decay}) was found by limiting the search interval for \textit{p} to $[1.5, 2]$ and setting the bisection inequality to $v(10)-\beta<0 \implies b = p \quad \textrm{else} \quad a = p$.\\
The one-node solution (Figure \ref{fig:onenode}) was found by limiting the search interval for \textit{p} to $[2, 2.5]$ and reversing the inequality sign. In addition, to avoid division by $0$ the spatial interval started at $10^{-6}$.

\subsection{Nonlinear Finite-Difference Method}
\begin{figure}[t]
    \centering
    \includegraphics[width=0.5\textwidth]{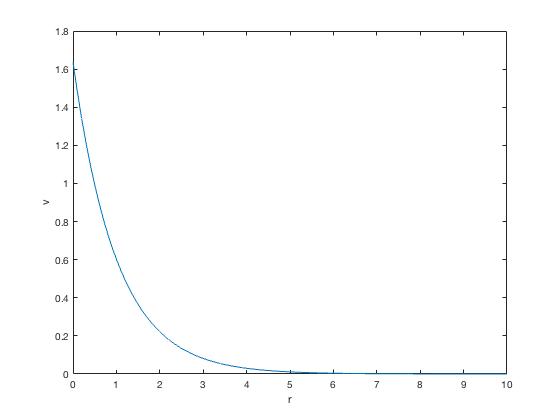}
    \caption{Initial $\vec{w}$ for Decaying Solution}
    \label{fig:decayinit}
\end{figure}
\begin{figure}[t]
    \centering
    \includegraphics[width=0.5\textwidth]{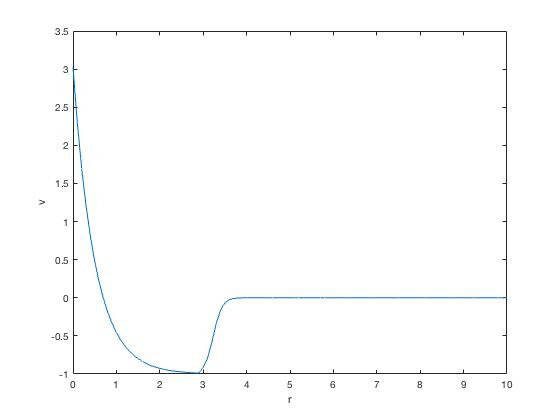}
    \caption{Initial $\vec{w}$ for One-Node Solution}
    \label{fig:onenodeinit}
\end{figure}
This method approximates the differential equation by replacing the derivatives in the equation with appropriate difference-quotient approximations. The interval $[0, 10]$ is subdivided into \textit{N} intervals with mesh points $a = x_1, \cdots , x_N, x_N+1 = b$ equally spaced by $\displaystyle h = \frac{b - a}{N}$. At each point $x_i$ the solution $v(r)$ is approximated by $w(x_i)$. \\
A system of N nonlinear equations must be solved:
\begin{equation}
\label{eq:3}
	\centering \displaystyle \frac{1}{h}(-\frac{3}{2}w_0+2w_1-\frac{1}{2}w_2)-0=0
\end{equation}
\begin{equation}
\label{eq:4}
	\centering \displaystyle -w_i+2w_{i+1}-w_{i+2}+h^2\Big[-\frac{1}{r}\frac{w_{i+2}-w_i}{2h}+w_i-2w_i^3\Big]=0
\end{equation}
for $i \in [1, N-1]$ where $w_{N+1}=0$.\\
Equation \ref{eq:3} is the Three-Point Endpoint formula which incorporates the boundary condition $\displaystyle \frac{dv(0)}{dr}=0$ from the BVP \ref{eq:1}. The remaining equations \ref{eq:4} are the centered-difference formulae of the Finite-Difference method.\\
For this method to work an initial guess for $\displaystyle \vec{w}$ had to be established. This initial guess had to be close enough to the solution to which the Finite-Difference method converges. For the decaying solution (Figure \ref{fig:decay}) an initial guess of
\begin{equation}
\label{eq:5}
	\centering \displaystyle  w(i) = 2e^{-0.1(i + 1)}
\end{equation}
was chosen which can be seen in Figure \ref{fig:decayinit}. For the one-node solution a piecewise function was used
\[
w(i)=
\begin{cases}
6e^{-0.2(i + 1)} - 1 &\text{if \ $(i-1)h) < 3$},\\[2ex]
(1 + e^{-i + \dfrac{N}{3}})^{-1}-1 &\text{if \ $(i-1)h) > 3$}.
\end{cases}
\]
which can be seen in Figure \ref{fig:onenodeinit}. The reason of the piecewise function was to establish a better initial guess that is closely related to the final approximation.\\
In order for the tridiagonal Jacobian matrix to be setup $\displaystyle \frac{\partial f}{\partial v}$ and $\displaystyle \frac{\partial f}{\partial v'}$ had to be computed where
\begin{equation}
\label{eq:5}
	\centering \displaystyle f(r, v, v') = v'' = -\frac{1}{r}\frac{dv}{dr}+v-2v^3
\end{equation}
\begin{equation}
\label{eq:5}
	\centering \displaystyle \frac{\partial f}{\partial v} = 1 - 6v^2
\end{equation}
\begin{equation}
\label{eq:5}
	\centering \displaystyle \frac{\partial f}{\partial v'} = -\frac{1}{r}
\end{equation}

\section{Results}
All experiments were carried out over a spatial interval of $[0, 10]$. For both the Shooting and Finite-difference method the amount of iterations required to reach the decaying and one-node solution are portrayed in table \ref{table:results1}, \ref{table:results2} and \ref{table:results3} (where the different tables represent varying error tolerances).

\begin{table}[t]
\centering
 \begin{tabular}{ c c c } 
 \hline
 & Decaying & One-Node \\ [0.5ex]
 \hline\hline
 Shooting & 29 & 25 \\
 Finite-Difference & 18 & 20 \\[0.5ex]
 \hline
 \end{tabular}
 \caption{Iterations of methods for the different solutions. Error $= 10^{-6}$}
\label{table:results1}
\end{table}

\begin{table}[t]
\centering
 \begin{tabular}{ c c c } 
 \hline
 & Decaying & One-Node \\ [0.5ex]
 \hline\hline
 Shooting & 36 & 35 \\
 Finite-Difference & 27 & 28 \\[0.5ex]
 \hline
 \end{tabular}
 \caption{Iterations of methods for the different solutions. Error $= 10^{-9}$}
\label{table:results2}
\end{table}

\begin{table}[t]
\centering
 \begin{tabular}{ c c c } 
 \hline
 & Decaying & One-Node \\ [0.5ex]
 \hline\hline
 Shooting & \textit{NA} & \textit{NA} \\
 Finite-Difference & 36 & 37 \\[0.5ex]
 \hline
 \end{tabular}
 \caption{Iterations of methods for the different solutions. Error $= 10^{-12}$}
\label{table:results3}
\end{table}

\section{Discussion}
From the results in can be seen that both the Shooting and Finite-Difference method require relatively few iterations to converge to an approximate solutions within the provided error tolerance. The Finite-Difference method outperforms the shooting method in regards to the amount of iterations required. In addition the Finite-Difference method has better stability characteristics as can be seen in Table \ref{table:results3} as it is able to find an approximate solution within an error tolerance of $10^{-12}$ which the shooting method is not able to do (i.e. it stops iterating once the max iterations count has been reached). \\
For the Shooting method it was important to centre the initial guess for the parameter $p$ to the defined intervals described in section 2.1 or else oscillating solutions were found. In addition solutions were still found over an extended spatial interval such as $[0, 30]$, however extending further resulted in instabilities.\\
The Finite-Difference method contains more hyper-parameters compared to the Shooting method, namely \textit{mesh spacing}, \textit{initial guess} and \textit{length of spatial interval}. An important observation made was that all these hyper-parameters are tightly linked; Changing one of them resulted in an oscillating solution. All the experiments used a mesh spacing of $100$ subintervals (i.e. $101$ mesh points). Increasing or decreasing this number resulted in undesirable solutions. The spatial length was kept at $[0, 10]$ as an increased spatial length suffered from the same issue. Most important was the observation of the sensitivity in regard to choosing an initial guess for $\vec{w}$. The decaying solution had an \textit{easier} initial guess compared to the one-node solution which required the construction of a piecewise function. The initial guess for both solutions were found by trial and error.

\section{Conclusion}
Both the Shooting and Finite-Difference method are powerful methods to obtain approximate solutions to nonlinear BVPs. The results show that the finite-difference method is numerically more stable and converges using less iterations compared to the Shooting method. However, the Finite-Difference method is more complicated to implement. Having more hyper-parameters than the Shooting method it is more difficult to fine tune to converge to the wanted solution as the mesh spacing, spatial interval and initial guess had to be determined via trial and error. It is to be concluded that the Finite-Difference method is to be used when approximate solutions with a low error threshold are to be found, otherwise the Shooting method is a better option.

\bibliographystyle{named}
\bibliography{ijcai17}

\section{Appendix}
\subsection{Shooting Method}
\begin{lstlisting}
% Decay sol in [1.5 2] @ 1.5609 when (v(end, 1) - beta < 0)
% One-node sol in [1.5 2] @ 2.3564 when (v(end, 1) - beta > 0)

function Shooting()
% Define the necessary variables to run the simulation
TOL = 10^-6; % Error threshold
MAX = 100;
END_POINT = 30; % End of interval
beta = 0; % limit as v -> inf
interval = [2 2.5]; % To search for y IC
IC = 0; % y' IC
[r, v] = shoot(interval, IC, END_POINT, beta, TOL, MAX);
disp('Completed');
print(r, v);
end

% Run the Nonlinear Shooting with Newton's Method iteratively
% Until the error falls below the defined Tollerance
function [r, v] = shoot(interval, IC, END_POINT, beta, TOL, MAX)
a = interval(1);
b = interval(2);
tk = (a+b)/2;
[r, v] = ode45(@ode, [10^-6, END_POINT], [tk, IC]);
iterations = 1;

while(abs(v(end, 1) - beta) > TOL)
    tk = (a + b) / 2;
    [r, v] = ode45(@ode, [10^-6, END_POINT], [tk, IC]);
    if(v(end, 1) - beta > 0)
        b = tk;
    else
        a = tk;
    end
    iterations = iterations + 1;
    if(iterations > MAX) 
        return;
    end
end
disp(iterations);
end

% Display solution
function print(t, y)
plot(t, y(:,1),'-')
xlabel('r');
ylabel('v');
end

% Original second order BVP expressed as a system of first-order equations
function dydt = ode(t, y)
dydt = [y(2); -(1/t)*y(2) + y(1) - 2*y(1)^3];
end	
\end{lstlisting}
\subsection{Finite-Difference Method}
\begin{lstlisting}
function Newton()
[r, v] = newton(100, 0, 10, 0, 0, -1, 10^-9);
disp('Completed');
print(r, v);
end

% f'(a) = alpha
% f(b) = beta
function [x, y] = newton(N, a, b, alpha, beta, M, TOL)
h = (b - a) / (N);
x = linspace(a, b, N + 1);
w = zeros;
% First solution
%for i = 1 : N
%    w(i) = exp(-(i + 1)*0.1)*2;
%end
% Second solution
for i = 1 : N
    if(a + ((i-1) * h) < 3)
        w(i) = exp(-(i + 1) * 0.2) * 6 - 1;
    else
        w(i) = 1 / (1 + exp(-i + N/3))-1;
    end  
end
w(N + 1) = beta;

k = 0;
while(k <= M)
    % Initialize column vector F
    F(1, 1) = (1/h) * (-(3/2) * w(1) + 2 * w(2) - (1/2) * w(3)) - alpha;
    for i = 1 : N - 1
        F(i + 1, 1) = - w(i) + 2 * w(i + 1) - w(i + 2) + (h ^ 2) * f(x(i + 1), w(i + 1), (w(i + 2) - w(i)) / (2 * h));
    end
    
    % Initialize tridiagonal Jacobian matrix
    J(1,1) = - 3 / (2 * h);
    J(1, 2) = 2 / h;
    J(1, 3) = - 1 / (2 * h);
    % Top diagonal
    for j = 3 : N
        i = j - 1;
        J(i, j) = - 1 + (h / 2) * fydif(x(i));
    end
    % Main diagonal
    for j = 2 : N
        i = j;
        J(i, j) = 2 + (h ^ 2) * fy(w(i));
    end
    % Bottom diagonal
    for j = 2 : N - 1
        i = j + 1;
        J(i, j) =  - 1 - (h / 2) * fydif(x(i));
    end
    
    % Solve for v
    v = inv(J) * - F;
    
    % Exit Loop if L2 Norm of v is less than the tolerance
    if(norm(v) < TOL)
        break;
    end
    
    % Update w based on adjustments from v
    for i = 1 : N
        w(i) = w(i) + v(i);
    end
    k = k + 1;
    
end
disp(k);
y = w(:);

end

function val = f(x, y, ydif)
val = - (1 / x) * ydif + y - 2 * (y ^ 3);
end

function val = fydif(x)
val = - (1 / x);
end

function val = fy(y)
val = 1 - 6 * y ^ 2;
end

% Display solution
function print(t, y)
plot(t, y(:,1),'-')
xlabel('r');
ylabel('v');
end	
\end{lstlisting}

\end{document}